\documentclass{article}
\usepackage{amsfonts, amssymb, amsmath}

\begin{document}

\begin{center}
{\Large The number of }$k${\large -potent elements in the quaternion algebra }$%
\mathbb{H}_{\mathbb{Z}_{p}}$%
\begin{equation*}
\end{equation*}

\textbf{Cristina FLAUT and Andreea BAIAS}%
\begin{equation*}
\end{equation*}
\end{center}

\textbf{Abstract.} {\small In this paper we count the number of }$k${\small %
-potent elements over }$H_{\mathbb{Z}_{p}}${\small , the quaternion algebra
over }$Z_{p}${\small , and we give a descriptive formula for the general
case. For }$k\in \{3,4,5\}${\small , we give an explicit formula for these
values. Moreover, as an application, we count the number of solutions of the
equation }$x^{k}=1${\small \ over }$H_{\mathbb{Z}_{p}}${\small .} 
\begin{equation*}
\end{equation*}

\textbf{2023-MSC: 17A45,17A75}

\begin{equation*}
\end{equation*}

\textbf{1.} \textbf{Introduction}%
\begin{equation*}
\end{equation*}

We consider the field $\mathbb{Z}_{p},p$ a prime odd number and the
quaternion algebra over $\mathbb{Z}_{p}=\{\widehat{0},\widehat{1},\widehat{2}%
,...,\widehat{p-1}\}$, denoted $\mathbb{H}_{\mathbb{Z}_{p}}=\left( \frac{%
-1,-1}{\mathbb{Z}_{p}}\right) $, with basis $\{\mathbf{1,i.j.k}\}$ and
multiplication given in the following table:%
\begin{equation*}
\begin{tabular}{l|llll}
$\cdot $ & $\mathbf{1}$ & $\mathbf{i}$ & $\mathbf{j}$ & $\mathbf{k}$ \\ 
\hline
$\mathbf{1}$ & $\mathbf{1}$ & $\mathbf{i}$ & $\mathbf{j}$ & $\mathbf{k}$ \\ 
$\mathbf{i}$ & $\mathbf{i}$ & $\mathbf{-1}$ & $\mathbf{k}$ & $\mathbf{-j}$
\\ 
$\mathbf{j}$ & $\mathbf{j}$ & $\mathbf{-k}$ & $\mathbf{-1}$ & $\mathbf{i}$
\\ 
$\mathbf{k}$ & $\mathbf{k}$ & $\mathbf{j}$ & $\mathbf{-i}$ & $\mathbf{-1}$%
\end{tabular}%
\end{equation*}

In this paper, we count the number of $k$-potent elements over $\mathbb{H}_{%
\mathbb{Z}_{p}}$ and we give a descriptive formula for the general case. For 
$k\in \{3,4,5\}$, we give an explicit formula for these values. The proposed
methods give us all elements with these properties. Moreover, as an
application, we count the number of solutions of the equation $x^{k}=1$ over 
$\mathbb{H}_{\mathbb{Z}_{p}}$. For this purpose, we used computer as a tool
to check and understand the behavior of these elements in each studied case.
This allowed us to give a correct mathematical proof for obtained relations.

In the following, in this section, we present some properties of quaternions
over $\mathbb{Z}_{p}$ useful for us in entire paper. For other details
regarding properties of quaternions over an arbitrary field, the reader is
referred to \textbf{[}Sc; 66\textbf{], [}Vo; 21\textbf{], [}FSF; 19\textbf{%
], }p. 431-449, etc.

If $q\in \mathbb{H}_{\mathbb{Z}_{p}},$ $%
q=q_{0}+q_{1}f_{1}+q_{2}f_{2}+q_{3}f_{3},$ then 
\begin{equation*}
\overline{q}=q_{0}-q_{1}f_{1}-q_{2}f_{2}-q_{3}f_{3}
\end{equation*}%
is called the \textit{conjugate} of the element $q.$ For $q\in \mathbb{H}_{%
\mathbb{Z}_{p}},$ we consider the following elements:

\begin{equation*}
\mathbf{t}\left( q\right) =q+\overline{q}\in \mathbb{H}_{\mathbb{Z}_{p}}
\end{equation*}%
and

\begin{equation*}
\,\mathbf{n}\left( q\right) =q\overline{q}%
=q_{0}^{2}+q_{1}^{2}+q_{2}^{2}+q_{3}^{2}\in \mathbb{H}_{\mathbb{Z}_{p}},
\end{equation*}%
called the \textit{trace}, respectively, the \textit{norm} of the element $%
q\in \mathbb{H}_{\mathbb{Z}_{p}}$. Shortly, the norm and the trace of the
element $q$ will be denoted $\mathbf{n}_{q}$ and $\mathbf{t}_{q}$.\thinspace
\thinspace It\thinspace \thinspace \thinspace follows\thinspace \thinspace
\thinspace that$\,\,$%
\begin{equation*}
\left( q+\overline{q}\right) q\,=q^{2}+\overline{q}q=q^{2}+\mathbf{n}\left(
q\right) \cdot 1
\end{equation*}%
and\thinspace \thinspace 
\begin{equation*}
q^{2}-\mathbf{t}\left( q\right) q+\mathbf{n}\left( q\right) =0,\forall q\in 
\mathbb{H}_{\mathbb{Z}_{p}},\,
\end{equation*}%
therefore the generalized quaternion algebra is a \textit{quadratic algebra}.

Since the field $\mathbb{Z}_{p}$ is a finite field, due to the Wedderburn's
Theorem, the quaternion algebra $\mathbb{H}_{\mathbb{Z}_{p}}$ is allways a
non division algebra or a split algebra. Also, we must remark that the
multiplicative group $\left( \mathbb{Z}_{p}^{\ast },\cdot \right) $ is a
cyclic group, for $p$ a prime number.\smallskip

\textbf{Definition 1.} An element $x$ in a ring $R$ is called \textit{%
nilpotent }if we can find a positive integer $n$ such that $x^{n}=0$. The
number $n$, the smallest with this property, is called the \textit{%
nilpotency index}. \ An element $x$ in a ring $R$ is called $k$\textit{%
-potent}, for $k>1$, a positive integer, if $k$ is the smallest number such
that $x^{k}=x$. The number $k$ is called the $k$\textit{-potency index}. For 
$k=2,$ we have idempotent elements, for $k=3$, we have tripotent elements,
and so on.\smallskip\ 

In the paper [MS; 11], authors counted the number of zero divisors in $%
\mathbb{H}_{\mathbb{Z}_{p}},$ that means the number of nonzero elements $%
x\in \mathbb{H}_{\mathbb{Z}_{p}}$ such that $\mathbf{n}\left( x\right) =0$.

In the following we will count the number of $k$-potent elements in
\smallskip $\mathbb{H}_{\mathbb{Z}_{p}}$ for $k\geq 3$. We will denote this
number with $\mathcal{I}_{p}\left( k\right) $.\smallskip

\textbf{Proposition 2.} i) ([MS; 11], Theorem 2.1) \textit{The number of
zero divisors in} $\mathbb{H}_{\mathbb{Z}_{p}}$ \textit{is} $p^{3}+p^{2}-p$.

ii) ([MS; 11], Theorem 2.1) \textit{The number of idempotent elements in} $%
\mathbb{H}_{\mathbb{Z}_{p}}$ \textit{is} $p^{2}+p+2$.\smallskip

We must remark that the above value for idempotents elements includes $0$
and $1$ as idempotents.\smallskip

\textbf{Remark 3. ([}FB; 24], Remark 2, \textit{iv}\textbf{))} In $\mathbb{H}%
_{\mathbb{Z}_{p}}$, a nilpotent element $x\in \mathbb{H}_{\mathbb{Z}_{p}}$
has $2$ as a nilpotency index and $\mathbf{n}\left( x\right) =\mathbf{t}%
\left( x\right) =0$.\smallskip \smallskip

In the paper [AD; 12], by using some results obtained in [FH; 58], authors
proved the following result.\smallskip

\textbf{Proposition 4.} ([AD; 12], Remark 2.4) \textit{The number of
nilpotent elements in} $\mathbb{H}_{\mathbb{Z}_{p}}$ \textit{is} $p^{2}$%
.\smallskip \smallskip

\textbf{Proposition 5.} \textit{The number of solutions in} $\mathbb{Z}_{p}$ 
\textit{of \ the equation}%
\begin{equation*}
x_{1}^{2}+x_{2}^{2}+x_{3}^{2}=0
\end{equation*}%
\textit{is} $p^{2}$.

\textbf{Proof.} The element $x=x_{1}\mathbf{i}+x_{2}\mathbf{j}+x_{3}\mathbf{k%
}$ has $\mathbf{t}_{x}=\mathbf{n}_{x}=0$, therefore it is nilpotent and we
apply the above proposition.\smallskip\ 

The above sentence also gives us a method to obtain all nilpotent elements
in $\mathbb{H}_{\mathbb{Z}_{p}}$, namely, all quaternions $x\in \mathbb{H}_{%
\mathbb{Z}_{p}}$ such that $\mathbf{t}_{x}=\mathbf{n}_{x}=0$.\smallskip\ 

\textbf{Remark 6.} From the above results, we have that, in $\mathbb{H}_{%
\mathbb{Z}_{p}}$, the number of elements with zero norm and nonzero trace is 
$\left( p^{3}+p^{2}-p\right) -p^{2}=p^{3}-p=p\left( p-1\right) \left(
p+1\right) $.\smallskip

Let $x,y\in \mathbb{H}_{\mathbb{Z}_{p}}$ be two quaternions, $x\neq y$, such
that $\mathbf{n}_{x}=\mathbf{n}_{y}=0$ and $\mathbf{t}_{x}\neq 0,\mathbf{t}%
_{y}\neq 0,\mathbf{t}_{x}\neq \mathbf{t}_{y}$. We denote with $A_{x}=\{w\in 
\mathbb{H}_{\mathbb{Z}_{p}},\mathbf{n}_{w}=0$ and$~$\ $\mathbf{t}_{w}=%
\mathbf{t}_{x}$ $\}$, $B_{y}=\{z\in \mathbb{H}_{\mathbb{Z}_{p}},\mathbf{n}%
_{z}=0$ and$~$\ $\mathbf{t}_{z}=\mathbf{t}_{y}$ $\}$.\smallskip\ 

\textbf{Proposition 7.} \textit{With the above notations, the sets} $A_{x}$ 
\textit{and} $B_{y}$ \textit{have the same number of elements, namely} 
\begin{equation*}
\left\vert A_{x}\right\vert =\left\vert B_{y}\right\vert =p\left( p+1\right)
.
\end{equation*}

\textbf{Proof.} Since $\mathbf{t}_{x}\neq 0,\mathbf{t}_{y}\neq 0$, we have $%
a,b\in \mathbb{H}_{\mathbb{Z}_{p}}$, invertible elements, such that $\mathbf{%
t}_{x}=a\mathbf{t}_{y}$ and $\mathbf{t}_{y}=b\mathbf{t}_{x}$. The map $%
\varphi :A_{x}\rightarrow B_{y},\varphi \left( v\right) =bv$ is an injective
map, then $\left\vert A_{x}\right\vert \leq \left\vert B_{y}\right\vert $.
Indeed, $bv$ has norm equal with zero. The map $\eta :B_{y}\rightarrow
A_{x},\eta \left( r\right) =ar$ is also an injective map, therefore $%
\left\vert B_{y}\right\vert \leq $ $\left\vert A_{x}\right\vert $. From
here, we have that $\left\vert A_{x}\right\vert =\left\vert B_{y}\right\vert 
$. From Remark 5, the number of elements with zero norm and nonzero trace is 
$p\left( p-1\right) \left( p+1\right) $ and the traces of these elements
have value in the set $\{\widehat{1},\widehat{2},...,\widehat{p-1}\}$. It
results that we have $p-1$ sets $A_{x}$, all having the same number of
elements, with $\mathbf{t}_{x}\in \{\widehat{1},\widehat{2},...,\widehat{p-1}%
\}$, therefore $\left\vert A_{x}\right\vert =\left\vert B_{y}\right\vert
=p\left( p+1\right) $.\smallskip

In [W; 22], the author proved the following results regarding the number of
solutions of particular types of equation over $\mathbb{Z}_{p}$.\smallskip

\textbf{Proposition 8. }([W; 22], Theorem 4.1) \textit{Let} $g\in \mathbb{Z}%
_{p}$, $g\neq 0$, \textit{then the equation} $x^{2}+y^{2}=g$ \textit{has }$%
p-\sin \frac{p\pi }{2}$ \textit{solutions in }$\mathbb{Z}_{p}$\textit{%
.\smallskip }

\textbf{Proposition 9.} \textit{Let} $g\in \mathbb{Z}_{p}$, $g\neq 0$, $g$ 
\textit{is not a square in} $\mathbb{Z}_{p}$. \textit{Then, the equation} 
\begin{equation}
x^{2}+y^{2}+z^{2}=g  \tag{1}
\end{equation}%
\textit{has} $p\left( p-\sin \frac{p\pi }{2}\right) $ \textit{%
solutions.\smallskip }

\textbf{Proof.} For each $z\in \{0,1,2,...,p-1\}$, equation $\left( 1\right) 
$ has $\left( p-\sin \frac{p\pi }{2}\right) $, from the above proposition,
therefore $p\left( p-\sin \frac{p\pi }{2}\right) $ is the number of all
solutions.\smallskip

\textbf{Proposition 10.} \textbf{\ }([W; 22], Theorem 4.2) \textit{Equation} 
$x^{2}+y^{2}=0$ \textit{has} $\mathcal{N}_{0}$ \textit{solutions, where} $%
\mathcal{N}_{0}=2p-1$\textit{, if} $p$ \textit{is a prime of the form} $%
4l+1,l\in \mathbb{Z}$ \textit{and }$\mathcal{N}_{0}=1$\textit{, if} $p$ 
\textit{is a prime of the form} $4l+3,l\in \mathbb{Z}$.\smallskip\ 

\textbf{Proposition 11.} \textit{Let} $g\in \mathbb{Z}_{p}$, $g\neq 0$, $g$%
\textit{\ is a square in} $\mathbb{Z}_{p}$. \textit{Then, the equation} 
\begin{equation*}
x^{2}+y^{2}+z^{2}=g
\end{equation*}%
\textit{has} $2\mathcal{N}_{0}+\left( p-2\right) \left( p-\sin \frac{p\pi }{2%
}\right) $ \textit{solutions.\smallskip }

\textbf{Proof.} We have $\frac{p-1}{2}$ nonzero perfect squares in $\mathbb{Z%
}_{p}$. Indeed, if we consider the following group morphism $\varphi :%
\mathbb{Z}_{p}^{\ast }\rightarrow (\mathbb{Z}_{p}^{\ast })^{2}$, $\varphi
\left( x\right) =x^{2}$, we have that $\varphi $ is surjective with \textit{%
Ker}$\varphi =\{1,p-1\}$. From Fundamental Theorem of Isomorphism we have $%
\mathbb{Z}_{p}^{\ast }/$\textit{Ker}$\varphi \simeq (\mathbb{Z}_{p}^{\ast
})^{2}$, therefore the cardinal of $(\mathbb{Z}_{p}^{\ast })^{2}$ is $\frac{%
p-1}{2}$. If $g=h^{2}$ and $z=h$ or $z=-h$, then equation $\left( 1\right) $
becomes $x^{2}+y^{2}=0$ and has $\mathcal{N}_{0}$ solutions. From here, the
number of solutions is $2\mathcal{N}_{0}+\left( p-2\right) \left( p-\sin 
\frac{p\pi }{2}\right) $.\smallskip

\textbf{Remark 12.} If $p$ is a prime odd number, then $\sin \frac{p\pi }{2}%
=1$, if $p=4l+1,l\in \mathbb{Z}$ and $\sin \frac{p\pi }{2}=-1$, if $%
p=4l+3,l\in \mathbb{Z}$.

\begin{equation*}
\end{equation*}

\textbf{2.} \textbf{The number of} $3$\textbf{-potent elements in quaternion
algebra} $\mathbb{H}_{\mathbb{Z}_{p}}$%
\begin{equation*}
\end{equation*}

An element $x\in $ $\mathbb{H}_{\mathbb{Z}_{p}}$ is called tripotent if $3$
is the smallest positive integer such that $x^{3}=x$. From here, we have two
cases: $\mathbf{n}_{x}=0$ or $\mathbf{n}_{x}\neq 0$ and $\mathbf{n}%
_{x}^{2}=1 $.

\textbf{Case 1,} $\mathbf{n}_{x}=0$. In this situation, we have that $%
\mathbf{t}_{x}\neq 0$, otherwise the element $x$ is nilpotent. \ From [FB;
24], Proposition 1, it results that $\mathbf{t}_{x}$ is $3$-potent in $%
\mathbb{Z}_{p}$, therefore $\mathbf{t}_{x}^{2}=1$ and $\mathbf{t}_{x}$ is an
element of order $2$ in the multiplicative group $\mathbb{Z}_{p}^{\ast }$.
Let $\Theta _{p}\left( 2\right) $ be the number of elements of order $2$ in $%
\mathbb{Z}_{p}^{\ast }$.

\textbf{Remark 13.} In $\mathbb{Z}_{p}$, we have only one element of order $%
2 $. Indeed, if $a\in \mathbb{Z}_{p}$ is of order $2$, it results that $%
p\mid \left( a^{2}-1\right) $, that means $p\mid \left( a-1\right) $ or $%
p\mid \left( a+1\right) $, with \thinspace $a\leq p-1$. From here, we have $%
a=1$ or $a=p-1$ and $\Theta _{p}\left( 2\right) =1$.\smallskip\ 

From Proposition 1, since for a fixed $\mathbf{t}_{x}$ and $\mathbf{n}_{x}=0$%
, we have $p\left( p+1\right) $ elements in $\mathbb{H}_{\mathbb{Z}_{p}}$
satifying relation 
\begin{equation*}
\mathbf{n}_{x}=x_{0}^{2}+x_{1}^{2}+x_{2}^{2}+x_{3}^{2}=0,
\end{equation*}
it results that the number of tripotent elements $\mathbb{H}_{\mathbb{Z}%
_{p}}~$with $\mathbf{n}_{x}=0,\mathbf{t}_{x}\neq 0\,$\ and $\mathbf{t}%
_{x}^{2}=1$ is $p\left( p+1\right) $.\smallskip\ 

\textbf{Case 2,} $\mathbf{n}_{x}\neq 0$, $\mathbf{n}_{x}^{2}=1$ and $x^{2}=1$%
. Since $x^{2}-\mathbf{t}_{x}x+\mathbf{n}_{x}=0$, it results $\mathbf{t}%
_{x}x=1+\mathbf{n}_{x}$.

i) If $\mathbf{t}_{x}\neq 0$, therefore $x=\mathbf{t}_{x}^{-1}\left( 1+%
\mathbf{n}_{x}\right) \in \mathbb{Z}_{p}$ is an element of order $2$ in $%
\mathbb{Z}_{p}$. From above remark, $x=1$ or $x=p-1$, then we have only one
tripotent element, $x=p-1$.

ii) If $\mathbf{t}_{x}=0$, then $1+\mathbf{n}_{x}=0$, therefore $\mathbf{n}%
_{x}=p-1$. Since $\mathbf{t}_{x}=0$, we have the element $x$ of the form $%
x=x_{1}\mathbf{i}+x_{2}\mathbf{j}+x_{3}\mathbf{k}$, such that $%
x_{1}^{2}+x_{2}^{2}+x_{3}^{2}=p-1$. We obtain the following relation in $%
\mathbb{Z}_{p}$.

\begin{equation}
1+x_{1}^{2}+x_{2}^{2}+x_{3}^{2}=0,\text{ }  \tag{2}
\end{equation}%
which has $p\left( p+1\right) $ number of solutions, due to the Proposition
6. Therefore, we have $p\left( p+1\right) $ tripotent elements in this case.
We get the following result.\smallskip

\textbf{Theorem 14.} \textit{The number of tripotent elements in} $\mathbb{H}%
_{\mathbb{Z}_{p}}$ \textit{is} 
\begin{equation}
\mathcal{I}_{p}\left( 3\right) =2p^{2}+2p+1.\smallskip  \tag{3}
\end{equation}

\textbf{Proof.} Indeed, by adding the obtained elements from Case 1 and Case
2, we have $\mathcal{I}_{p}\left( 3\right) =p\left( p+1\right) +p\left(
p+1\right) +1=2p^{2}+2p+1$.\smallskip

We remark that from the above proposed method, we can obtain all tripotent
elements as the solutions of \ equations given in Case 1 and Case 2. 
\begin{equation*}
\end{equation*}

\textbf{3.} \textbf{The number of} $4$\textbf{-potent elements in quaternion
algebra} $\mathbb{H}_{\mathbb{Z}_{p}}$

\begin{equation*}
\end{equation*}

An element $x\in $ $\mathbb{H}_{\mathbb{Z}_{p}}$ is called $4$-potent if $4$
is the smallest positive integer such that $x^{4}=x$. From here, we have two
cases: $\mathbf{n}_{x}=0$ or $\mathbf{n}_{x}\neq 0$ and $\mathbf{n}%
_{x}^{3}=1 $.

\textbf{Case 1,} $\mathbf{n}_{x}=0$. In this situation, we have that $%
\mathbf{t}_{x}\neq 0$, otherwise the element $x$ is nilpotent. From [FB;
24], Proposition 1, it results that $\mathbf{t}_{x}$ is $4$-potent in $%
\mathbb{Z}_{p}$, therefore $\mathbf{t}_{x}^{3}=1$ and $\mathbf{t}_{x}$ is an
element of order $3$ in the multiplicative group $\mathbb{Z}_{p}^{\ast }$.
Let $\Theta _{p}\left( 4\right) $ be the number of elements of order $3$ in
the multiplicative group $\mathbb{Z}_{p}^{\ast }$.

From Proposition 1, since for a fixed $\mathbf{t}_{x}$ and $\mathbf{n}_{x}=0$%
, we have $p\left( p+1\right) $ elements in $\mathbb{H}_{\mathbb{Z}_{p}}$
satifying relation 
\begin{equation*}
\mathbf{n}_{x}=x_{0}^{2}+x_{1}^{2}+x_{2}^{2}+x_{3}^{2}=0,
\end{equation*}
it results that the number of $4$-potent elements $\mathbb{H}_{\mathbb{Z}%
_{p}}~$with $\mathbf{n}_{x}=0,\mathbf{t}_{x}\neq 0\,$\ and $\mathbf{t}%
_{x}^{3}=1$ is $\Theta _{p}\left( 4\right) p\left( p+1\right) $.

\textbf{Case 2,} $\mathbf{n}_{x}\neq 0$, $\mathbf{n}_{x}^{3}=1$ and $x^{3}=1$%
.

i) $\mathbf{t}_{x}=0$. It results $x^{2}+\mathbf{n}_{x}=0$, therefore $%
x^{2}=-\mathbf{n}_{x}$. We obtain $-\mathbf{n}_{x}x=1$. Since $x^{3}=1$ and $%
x=-\mathbf{n}_{x}^{-1}$, we have $\mathbf{n}_{x}^{3}=-1$, false, since $%
\mathbf{n}_{x}^{3}=1$ and $p>2$. Therefore $\mathbf{t}_{x}\neq 0$.

ii) $\mathbf{t}_{x}\neq 0$. From $x^{3}=1$ and since $x^{2}-\mathbf{t}_{x}x+%
\mathbf{n}_{x}=0$, it results $1-\mathbf{t}_{x}x^{2}+\mathbf{n}_{x}x=0$ and
we obtain $\mathbf{t}_{x}x^{2}-\mathbf{n}_{x}x-1=0$. Since $\mathbf{t}%
_{x}\neq 0$, we have $x^{2}-\mathbf{t}_{x}^{-1}\mathbf{n}_{x}x-\mathbf{t}%
_{x}^{-1}=0$ and we obtain $\left( \mathbf{t}_{x}^{-1}\mathbf{n}_{x}-\mathbf{%
t}_{x}\right) x=-(\mathbf{n}_{x}+\mathbf{t}_{x}^{-1})$.

-If $\mathbf{t}_{x}^{-1}\mathbf{n}_{x}-\mathbf{t}_{x}\neq 0$, it results
that $x\in \mathbb{Z}_{p}$. Since it is an element of order $3$ in $\mathbb{Z%
}_{p}$, different from $1$, we obtain $\Theta _{p}\left( 4\right) \,$\ as
the number of $4$-potent elements in $\mathbb{H}_{\mathbb{Z}_{p}}$, in this
case.

-If $\mathbf{t}_{x}^{-1}\mathbf{n}_{x}=\mathbf{t}_{x}$, also $(\mathbf{n}%
_{x}+\mathbf{t}_{x}^{-1})=0$ and it results that $\mathbf{n}_{x}=-\mathbf{t}%
_{x}^{-1}$, therefore 
\begin{equation*}
\mathbf{n}_{x}=\mathbf{t}_{x}^{2}\text{, }\mathbf{n}_{x}\mathbf{t}_{x}=-1%
\text{ and }\mathbf{t}_{x}^{3}=-1\text{. }
\end{equation*}%
Let $\Upsilon _{p}\left( 4\right) $ be the number of solutions of the
equation $\mathbf{t}_{x}^{3}=-1$ in $\mathbb{Z}_{p}$. For $x=x_{0}+x_{1}%
\mathbf{i}+x_{2}\mathbf{j}+x_{3}\mathbf{k}$, we have $x_{0}=2^{-1}\mathbf{t}%
_{x}$ and relation $\mathbf{n}_{x}=\mathbf{t}_{x}^{2}$ implies 
\begin{equation*}
\left( 2^{-1}\mathbf{t}_{x}\right) ^{2}+x_{1}^{2}+x_{2}^{2}+x_{3}^{2}=%
\mathbf{t}_{x}^{2}.
\end{equation*}

We obtain 
\begin{equation*}
\mathbf{t}_{x}^{2}+(2x_{1})^{2}+(2x_{2})^{2}+(2x_{3})^{2}=4\mathbf{t}%
_{x}^{2},
\end{equation*}%
that means 
\begin{equation}
(2x_{1})^{2}+(2x_{2})^{2}+(2x_{3})^{2}=3\mathbf{t}_{x}^{2}.  \tag{4}
\end{equation}

\textbf{Remark 15.} Let $a,b\in \mathbb{Z}_{p}$, $a\neq b$. Therefore $%
3a^{2} $ and $3b^{2}$ are in the same time quadratic resudues in $\mathbb{Z}%
_{p}$ or not. Indeed, if $3a^{2}=\beta ^{2}$, we have that $3$ is a perfect
square in $\mathbb{Z}_{p}$ and $3b^{2}$ is also a perfect square and
vice-versa.\smallskip\ From here and from Proposition 9 and Proposition 11,
we have that the number of solutions of the equation $\left( 4\right) $ is
the same for the each given value of $\mathbf{t}_{x}$. We denote this number
with $\mathfrak{N}_{4}$. By using Propositions 9, 10, 11, we obtain the
following values for $\mathfrak{N}_{4}$:

-If $3\mathbf{t}_{x}^{2}$ is not a perfect square in$\mathbb{\ Z}_{p}$, then 
$\mathfrak{N}_{4}=$ $p\left( p-\sin \frac{p\pi }{2}\right) $.

-If $3\mathbf{t}_{x}^{2}$ is a perfect square in$\mathbb{\ Z}_{p}$, then $%
\mathfrak{N}_{4}=2\mathcal{N}_{0}+\left( p-2\right) \left( p-\sin \frac{p\pi 
}{2}\right) $, where $\mathcal{N}_{0}=2p-1$\textit{, }if $p$ \textit{i}s a
prime of the form $4l+1,l\in \mathbb{Z}$ and\textit{\ }$\mathcal{N}_{0}=1$%
\textit{, }if $p$ is a prime of the form $4l+3,l\in \mathbb{Z}$.\smallskip

\textbf{Proposition 16.} \textit{With the above notations, we have that the
number of} $4$\textit{-potent} \textit{elements in} $\mathbb{H}_{\mathbb{Z}%
_{p}}$ \textit{with} $\mathbf{n}_{x}=\mathbf{t}_{x}^{2}$, $\mathbf{n}_{x}%
\mathbf{t}_{x}=-1$ \textit{and} $\mathbf{t}_{x}^{3}=-1$ \textit{is} $%
\Upsilon _{p}\left( 4\right) \mathfrak{N}_{4}$.\smallskip

\textbf{Theorem 17.} \textit{The number of }$4$-\textit{potent elements in} $%
\mathbb{H}_{\mathbb{Z}_{p}}$ \textit{is} 
\begin{equation}
\mathcal{I}_{p}\left( 4\right) =\Theta _{p}\left( 4\right) \left(
p^{2}+p+1\right) +\Upsilon _{p}\left( 4\right) \mathfrak{N}_{4}.\smallskip 
\tag{5}
\end{equation}

\textbf{Proof.} Indeed, by adding the obtained elements from Case 1 and Case
2, we have $\mathcal{I}_{p}\left( 4\right) =\Theta _{p}\left( 4\right)
p\left( p+1\right) +\Theta _{p}\left( 4\right) +\Upsilon _{p}\left( 4\right) 
\mathfrak{N}_{4}$.\smallskip

\textbf{Example 18.} i) If $p=3$, there are not elements of order $3$ in $%
\mathbb{Z}_{3}^{\ast }$. \ Therefore, $\Theta _{p}\left( 4\right) =0$. We
also have $\Upsilon _{3}\left( 4\right) =1$. Equation $\left( 4\right) $
becomes 
\begin{equation*}
x_{1}^{2}+x_{2}^{2}+x_{3}^{2}=0\text{.}
\end{equation*}%
This equation represents all nonzero nilpotent elements, then has $%
3^{2}-1=2^{3}$ solutions in $\mathbb{Z}_{3}$:$~\left\{ 1,1,1\right\} $, $%
\left\{ 2,1,1\right\} $, $\left\{ 1,2,1\right\} $, $\left\{ 1,1,2\right\} $, 
$\left\{ 2,2,1\right\} $, $\left\{ 2,1,2\right\} $, $\left\{ 1,2,2\right\} $%
, $\left\{ 2,2,2\right\} $. Therefore, the number of $4$-potent elements in $%
\mathbb{H}_{\mathbb{Z}_{3}}$ is $\mathcal{I}_{3}\left( 4\right) =8$.

ii) If $p=5$, there are not elements of order $3$ in $\mathbb{Z}_{5}^{\ast 
\text{ }}$. Therefore, $\Theta _{5}\left( 4\right) =0$ and we have $\Upsilon
_{5}\left( 4\right) =1$, since the solution of equation $\mathbf{t}%
_{x}^{3}=-1$ in $\mathbb{Z}_{5}$ is $\left\{ 4\right\} $. Equation $\left(
4\right) $ becomes $-3\mathbf{t}%
_{x}^{2}+(2x_{1})^{2}+(2x_{2})^{2}+(2x_{3})^{2}=0$. It results, $-3\mathbf{t}%
_{x}^{2}+x_{1}^{2}+x_{2}^{2}+x_{3}^{2}=-3+x_{1}^{2}+x_{2}^{2}+x_{3}^{2}=0$,
since $\mathbf{t}_{x}=4$. We must find the number of solutions of equation 
\begin{equation*}
x_{1}^{2}+x_{2}^{2}+x_{3}^{2}=3.
\end{equation*}%
Since $3$ is not a square in $\mathbb{Z}_{5}$, from Proposition 9, we have
that the number of the above equation is $p\left( p-\sin \frac{p\pi }{2}%
\right) =5\left( 5-\sin \frac{5\pi }{2}\right) =5\ast 4=20.$ Therefore, the
number of $4$-potent elements in $\mathbb{H}_{\mathbb{Z}_{5}}$ is $\mathcal{I%
}_{5}\left( 4\right) =20$.

iii) If $p=7$, then $2$ and $4$ are elements of order $3$ in $\mathbb{Z}%
_{7}^{\ast }$. It results that $\Theta _{7}\left( 4\right) =2$. Moreover, we
have $\Upsilon _{7}\left( 4\right) =3$, with the set $\left\{ 3,5,6\right\} $%
. Equation $\left( 4\right) $ becomes

\begin{equation*}
4\mathbf{t}_{x}^{2}+(2x_{1})^{2}+(2x_{2})^{2}+(2x_{3})^{2}=0,
\end{equation*}%
that means 
\begin{equation*}
\mathbf{t}_{x}^{2}+x_{1}^{2}+x_{2}^{2}+x_{3}^{2}=0.
\end{equation*}%
This equation has $\mathfrak{N}_{4}=p\left( p+1\right) =56$ solutions for
each value of $\mathbf{t}_{x}$. We obtain $\mathcal{I}_{7}\left( 4\right)
=\Theta _{p}\left( 4\right) \left( p^{2}+p+1\right) +\Upsilon _{p}\left(
4\right) \mathfrak{N}_{4}=2\ast 57+3\ast 56=282$, which represents the
number of $4$-potent elements in $\mathbb{H}_{\mathbb{Z}_{7}}$.

iv) If $p=11$, then there are not elements of order $3$ in $\mathbb{Z}%
_{11}^{\ast }$ and $\Theta _{11}\left( 4\right) =0$. We have $\Upsilon
_{11}\left( 4\right) =1$, since $10^{3}=10=-1$. Equation $%
(2x_{1})^{2}+(2x_{2})^{2}+(2x_{3})^{2}=3\mathbf{t}_{x}^{2}$ becomes $%
(2x_{1})^{2}+(2x_{2})^{2}+(2x_{3})^{2}=3$. Since $4^{-1}=3$ mod $11$, we
obtain the equation $x_{1}^{2}+x_{2}^{2}+x_{3}^{2}=9$. From Proposition 10
and Proposition 11, we have that the number of $4$-potent elements in $%
\mathbb{H}_{\mathbb{Z}_{11}}$ is $\mathcal{I}_{11}\left( 4\right) =2+\left(
p-2\right) \left( p+1\right) =p\left( p-1\right) =110$, for $p=11$.

v) If $p=13$, then we have $3$ and $9$ elements of order $3$ in $\mathbb{Z}%
_{13}^{\ast }$ and $\Theta _{13}\left( 4\right) =2$. We have $\Upsilon
_{13}\left( 3\right) =3$, with the set $\{4,10,12\}$.

For the first part, we have $\Theta _{13}\left( 4\right) \left(
p^{2}+p+1\right) =2\ast \left( 13^{2}+13+1\right) =366$ elements$.$

For the second part, we have the following values.

For $\mathbf{t}_{x}^{2}=3$, equation $%
(2x_{1})^{2}+(2x_{2})^{2}+(2x_{3})^{2}=3\mathbf{t}_{x}^{2}$ becomes $%
(2x_{1})^{2}+(2x_{2})^{2}+(2x_{3})^{2}=9$. Since $4^{-1}=10$ mod $13$, we
obtain the equation $x_{1}^{2}+x_{2}^{2}+x_{3}^{2}=12$. Since $12$ is a
square modulo $13$, $5^{2}=12$ mod $13$, from Proposition 9, we have that
the number of the above equation is $2\mathcal{N}_{0}+\left( p-2\right)
\left( p-\sin \frac{p\pi }{2}\right) =2\ast (2p-1)+\left( p-2\right) \left(
p-\sin \frac{p\pi }{2}\right) =2\ast 25+11\ast 12=182.$

For $\mathbf{t}_{x}^{2}=9$, equation $%
(2x_{1})^{2}+(2x_{2})^{2}+(2x_{3})^{2}=3\mathbf{t}_{x}^{2}$ becomes $%
(2x_{1})^{2}+(2x_{2})^{2}+(2x_{3})^{2}=9$. Since $4^{-1}=10$ mod $13$, we
obtain the equation $x_{1}^{2}+x_{2}^{2}+x_{3}^{2}=12$. Since $12$ is a
square modulo $13$, $5^{2}=12$ mod $13$, From Proposition 10 and Proposition
11, we have that the number of the above equation is $2\mathcal{N}%
_{0}+\left( p-2\right) \left( p-\sin \frac{p\pi }{2}\right) =$\newline
$=2\ast (2p-1)+\left( p-2\right) \left( p-\sin \frac{p\pi }{2}\right) =2\ast
25+11\ast 12=182$, the same value.

For $\mathbf{t}_{x}^{2}=1$, equation $%
(2x_{1})^{2}+(2x_{2})^{2}+(2x_{3})^{2}=3\mathbf{t}_{x}^{2}$ becomes $%
(2x_{1})^{2}+(2x_{2})^{2}+(2x_{3})^{2}=3$. Since $4^{-1}=10$ mod $13$, we
obtain the equation $x_{1}^{2}+x_{2}^{2}+x_{3}^{2}=4$. Since $4$ is a square
modulo $13$, $2^{2}=4$ mod $13$, From Proposition 10 and Proposition 11, we
have that the number of the above equation is $2\mathcal{N}_{0}+\left(
p-1\right) \left( p-\sin \frac{p\pi }{2}\right) =$\newline
$=2\ast (2p-1)+\left( p-2\right) \left( p-\sin \frac{p\pi }{2}\right) =2\ast
25+11\ast 12=182$, the same value.

Therefore, the number of $4$-potent elements in $\mathbb{H}_{\mathbb{Z}%
_{13}} $ is $\mathcal{I}_{13}\left( 4\right) =$ $912$.

vi) If $p=17$, then there are not elements of order $3$ in $\mathbb{Z}%
_{17}^{\ast }$ and $\Theta _{17}\left( 4\right) =0$. We have $\Upsilon
_{17}\left( 4\right) =1$, since $16^{3}=16=-1$. Equation $%
(2x_{1})^{2}+(2x_{2})^{2}+(2x_{3})^{2}=3\mathbf{t}_{x}^{2}$ becomes $%
(2x_{1})^{2}+(2x_{2})^{2}+(2x_{3})^{2}=3$. Since $4^{-1}=13$ mod $17$, we
obtain the equation $x_{1}^{2}+x_{2}^{2}+x_{3}^{2}=5$. Since $5$ is not a
square modulo $17$, from Proposition 9, we have that the number of the above
equation is $p\left( p-\sin \frac{p\pi }{2}\right) =17\left( 17-1\right)
=272 $, therefore the number of $4$-potent elements in $\mathbb{H}_{\mathbb{Z%
}_{17}}$ is $\mathcal{I}_{17}\left( 4\right) =272$.

vii) In the same way, for $p=19$, the number of $4$-potent elements in $%
\mathbb{H}_{\mathbb{Z}_{19}}$ is $\mathcal{I}_{19}\left( 4\right) =1902$.
For $p=23$, the number of $4$-potent elements in $\mathbb{H}_{\mathbb{Z}%
_{23}}$ is $\mathcal{I}_{23}\left( 4\right) =506$. For $p=29$, the number of 
$4$-potent elements in $\mathbb{H}_{\mathbb{Z}_{29}}$ is $\mathcal{I}%
_{29}\left( 4\right) =812$, etc.\smallskip

We remark that from the above proposed method, we can obtain all $4$-potent
elements as the solutions of \ equations given in Case 1 and Case 2. 
\begin{equation*}
\end{equation*}

\bigskip \textbf{4.} \textbf{The number of} $5$\textbf{-potent elements in
quaternion algebra} $\mathbb{H}_{\mathbb{Z}_{p}}$

\begin{equation*}
\end{equation*}

An element $x\in $ $\mathbb{H}_{\mathbb{Z}_{p}}$ is called $5$-potent if $5$
is the smallest positive integer such that $x^{5}=x$. From here, we have two
cases: $\mathbf{n}_{x}=0$ or $\mathbf{n}_{x}\neq 0$ and $\mathbf{n}%
_{x}^{4}=1 $.

\textbf{Case 1,} $\mathbf{n}_{x}=0$. In this situation, we have that $%
\mathbf{t}_{x}\neq 0$, otherwise the element $x$ is nilpotent. From [FB;
24], Proposition 1, it results that $\mathbf{t}_{x}$ is $5$-potent in $%
\mathbb{Z}_{p}$, therefore $\mathbf{t}_{x}^{4}=1$ and $\mathbf{t}_{x}$ is an
element of order $4$ in the multiplicative group $\mathbb{Z}_{p}^{\ast }$.
Let $\Theta _{p}\left( 5\right) $ be the number of elements of order $4$ in $%
\mathbb{Z}_{p}^{\ast }$.

From Proposition 1, since for a fixed $\mathbf{t}_{x}$ and $\mathbf{n}_{x}=0$%
, we have $p\left( p+1\right) $ elements in $\mathbb{H}_{\mathbb{Z}_{p}}$
satifying relation 
\begin{equation*}
\mathbf{n}_{x}=x_{0}^{2}+x_{1}^{2}+x_{2}^{2}+x_{3}^{2}=0,
\end{equation*}
it results that the number of $5$-potent elements in $\mathbb{H}_{\mathbb{Z}%
_{p}}~$with $\mathbf{n}_{x}=0,\mathbf{t}_{x}\neq 0\,$\ and $\mathbf{t}%
_{x}^{4}=1$ is $\Theta _{p}\left( 5\right) p\left( p+1\right) $.

\textbf{Case 2,} $\mathbf{n}_{x}\neq 0$, $\mathbf{n}_{x}^{4}=1$ and $x^{4}=1$%
.

i) If $\mathbf{t}_{x}=0$, we have $x^{2}=-\mathbf{n}_{x}$ and $\mathbf{n}%
_{x}x^{2}+1=0$, that means $x^{2}=-\mathbf{n}_{x}^{-1}$. It results $\mathbf{%
n}_{x}^{2}=1$. From here, we have that $p\mid \left( \mathbf{n}_{x}-1\right)
\left( \mathbf{n}_{x}+1\right) $. Therefore, \thinspace $\mathbf{n}_{x}=1$
or $\mathbf{n}_{x}=-1$ in $\mathbb{Z}_{p}$. If $\mathbf{n}_{x}=-1$, it
results $x^{2}-1=0$, that means $x^{2}=1$, which implies $x^{3}=x$, a
contradiction to the fact that $x$ is $5$-potent. From here, we obtain only
solution $\mathbf{n}_{x}=1$, that means $x^{2}=-1$. The number of $5$-potent
elements in this case is given by the number of solutions of equation%
\begin{equation*}
x_{1}^{2}+x_{2}^{2}+x_{3}^{2}=1.
\end{equation*}%
Therefore we have $2\mathcal{N}_{0}+\left( p-2\right) \left( p-\sin \frac{%
p\pi }{2}\right) $ elements, where $\mathcal{N}_{0}=2p-1$\textit{, }if $p$ 
\textit{i}s a prime of the form $4l+1,l\in \mathbb{Z}$ and\textit{\ }$%
\mathcal{N}_{0}=1$\textit{, }if $p$ is a prime of the form $4l+3,l\in 
\mathbb{Z}$.\smallskip

ii) If $\mathbf{t}_{x}\neq 0$, since $x^{2}-\mathbf{t}_{x}x+\mathbf{n}_{x}=0$%
, it results $1-\mathbf{t}_{x}x^{3}+\mathbf{n}_{x}x^{2}=0$ and we obtain $%
\mathbf{t}_{x}x^{3}-\mathbf{n}_{x}x^{2}-1=0$. From theses relations, we have 
$\mathbf{t}_{x}x\left( \mathbf{t}_{x}x-\mathbf{n}_{x}\right) -\mathbf{n}%
_{x}x^{2}-1=0$ and $(\mathbf{t}_{x}^{2}-\mathbf{n}_{x})x^{2}-\mathbf{t}_{x}%
\mathbf{n}_{x}x-1=0$.

-If $\mathbf{t}_{x}^{2}-\mathbf{n}_{x}=0$, it results $\mathbf{t}_{x}\mathbf{%
n}_{x}x+1=0$ and $x=-\left( \mathbf{t}_{x}\mathbf{n}_{x}\right) ^{-1}$,
therefore $x\in \mathbb{Z}_{p}$ is an element of order four, different from $%
1$. We obtain that $\mathbf{t}_{x}^{4}=1$ and, since $\mathbf{t}_{x}^{2}=%
\mathbf{n}_{x}$, we have $\mathbf{n}_{x}^{2}=1$. Then, as from above, we
have \thinspace $\mathbf{n}_{x}=1$ or $\mathbf{n}_{x}=-1$ in $\mathbb{Z}_{p}$%
. If $\mathbf{n}_{x}=1$, we have $\mathbf{t}_{x}^{2}=1,$ then $\mathbf{t}%
_{x}=1$ or $\mathbf{t}_{x}=-1$. If$~\mathbf{t}_{x}=\mathbf{n}_{x}=1$, we
obtain $x=-1$, therefore $x^{2}=1\,$, false, since $x$ is $5$-potent. If $%
\mathbf{t}_{x}=-1$ and $\mathbf{n}_{x}=1$, we have $x=1$, false. If $\mathbf{%
n}_{x}=-1$, then $\mathbf{t}_{x}^{2}=-1$. From $\mathbf{t}_{x}\mathbf{n}%
_{x}x+1=0$, we have $\mathbf{t}_{x}x=1$ and relation $x^{2}-\mathbf{t}_{x}x+%
\mathbf{n}_{x}=0$ becomes $x^{2}-1-1=0$, that means $x^{2}=2$. Since $%
x^{4}=4=1$ \textit{mod} $p$, we have $4=1$ \textit{mod} $p$, therefore $p=3$%
. Since in $\mathbb{Z}_{3}$ we have $\Theta _{3}\left( 5\right) =0\,$, there
are no solutions in this case.

-If $\mathbf{t}_{x}^{2}-\mathbf{n}_{x}\neq 0$, from the system $\left\{ 
\begin{array}{c}
x^{2}-\mathbf{t}_{x}x+\mathbf{n}_{x}=0 \\ 
(\mathbf{t}_{x}^{2}-\mathbf{n}_{x})x^{2}-\mathbf{t}_{x}\mathbf{n}_{x}x-1=0%
\end{array}%
\right. $, we have $(\mathbf{t}_{x}(\mathbf{t}_{x}^{2}-\mathbf{n}_{x})-%
\mathbf{t}_{x}\mathbf{n}_{x})x=1+\mathbf{n}_{x}(\mathbf{t}_{x}^{2}-\mathbf{n}%
_{x})$, then 
\begin{equation*}
\left( \mathbf{t}_{x}^{3}-2\mathbf{t}_{x}\mathbf{n}_{x}\right) x=\mathbf{n}%
_{x}\mathbf{t}_{x}^{2}-\mathbf{n}_{x}^{2}+1.
\end{equation*}

If $\mathbf{t}_{x}^{3}-2\mathbf{t}_{x}\mathbf{n}_{x}\neq 0$, we get $x\in 
\mathbb{Z}_{p}$. Since it is an element of order $4$ in $\mathbb{Z}_{p}$, we
obtain $\Theta _{p}\left( 5\right) \,$\ as the number of $4$-potent elements
in $\mathbb{H}_{\mathbb{Z}_{p}}$, in this case.

If $\mathbf{t}_{x}^{3}-2\mathbf{t}_{x}\mathbf{n}_{x}=0$, since $\mathbf{t}%
_{x}\neq 0$, we have $\mathbf{t}_{x}^{2}-2\mathbf{n}_{x}=0$. Also, we must
have $\mathbf{n}_{x}\mathbf{t}_{x}^{2}-\mathbf{n}_{x}^{2}+1=0$. From the
last two relations, we obtain $\mathbf{n}_{x}=2^{-1}\mathbf{t}_{x}^{2}$ \
and $\mathbf{n}_{x}^{2}=-1$, therefore $\mathbf{t}_{x}^{4}=-4$.

Let $\Upsilon _{p}\left( 5\right) $ be the number of solutions of the
equation $\mathbf{t}_{x}^{4}=-4$ in $\mathbb{Z}_{p}$. For $x=x_{0}+x_{1}%
\mathbf{i}+x_{2}\mathbf{j}+x_{3}\mathbf{k}$, we have $x_{0}=2^{-1}\mathbf{t}%
_{x}$ and relation $\mathbf{n}_{x}=2^{-1}\mathbf{t}_{x}^{2}$ implies 
\begin{equation*}
\left( 2^{-1}\mathbf{t}_{x}\right) ^{2}+x_{1}^{2}+x_{2}^{2}+x_{3}^{2}=2^{-1}%
\mathbf{t}_{x}^{2}.
\end{equation*}%
We obtain 
\begin{equation}
(2x_{1})^{2}+(2x_{2})+(2x_{3})=\mathbf{t}_{x}^{2}\text{.}  \tag{6}
\end{equation}

Since $\mathbf{t}_{x}^{2}$ is a perfect square in $\mathbb{Z}_{p}$, we
obtain the same number of solutions of equation $\left( 6\right) $, for each
value of $\mathbf{t}_{x}$. We denote this number with $\mathfrak{N}_{5}$. By
using Propositions 10 and 11, we obtain that $\mathfrak{N}_{5}=2\mathcal{N}%
_{0}+\left( p-2\right) \left( p-\sin \frac{p\pi }{2}\right) $, where $%
\mathcal{N}_{0}=2p-1$\textit{, }if $p$ \textit{i}s a prime of the form $%
4l+1,l\in \mathbb{Z}$ and\textit{\ }$\mathcal{N}_{0}=1$\textit{, }if $p$ is
a prime of the form $4l+3,l\in \mathbb{Z}$.\smallskip

\smallskip \textbf{Theorem 19.} \textit{The number of }$5$-\textit{potent
elements in} $\mathbb{H}_{\mathbb{Z}_{p}}$ \textit{is} 
\begin{equation}
\mathcal{I}_{p}\left( 5\right) =\Theta _{p}\left( 5\right) \left(
p^{2}+p+1\right) +(\Upsilon _{p}\left( 5\right) +1)\mathfrak{N}%
_{5}.\smallskip  \tag{7}
\end{equation}

\textbf{Proof.} Indeed, from the above, the number is $\mathcal{I}_{p}\left(
5\right) =\Theta _{p}\left( 5\right) p\left( p+1\right) +2\mathcal{N}%
_{0}+\left( p-2\right) \left( p-\sin \frac{p\pi }{2}\right) +\Theta
_{p}\left( 5\right) +\Upsilon _{p}\left( 5\right) \mathfrak{N}_{5}=$\newline
$=\Theta _{p}\left( 5\right) \left( p^{2}+p+1\right) +(\Upsilon _{p}\left(
5\right) +1)\mathfrak{N}_{5}$.\smallskip

\textbf{Example 20.} i) If $p=3$, then there are not elements of order $4$
in $\mathbb{Z}_{3}^{\ast }$. It results $\Theta _{3}\left( 5\right) =0$.
Since $2^{2}=2^{4}=1$, we also have $\Upsilon _{3}\left( 5\right) =0$. We
obtain $\mathcal{I}_{3}\left( 5\right) =\mathfrak{N}_{5}=2+1\ast 4=6$.

ii) If $p=5$, we have $\mathbf{2}^{4}\mathbf{=1},\mathbf{3}^{4}\mathbf{=1}%
,4^{2}=4^{4}=1$. Therefore $\Theta _{5}\left( 5\right) =2$. Since $-4=1$ in $%
\mathbb{Z}_{5}$, we have $\Upsilon _{5}\left( 5\right) =4$. We obtain $%
\mathcal{I}_{5}\left( 5\right) =2\ast \left( 25+5+1\right) +5\ast \left(
2\ast 9+3\ast 4\right) =212$.

iii) If $p=7$, we have $-4=3$ and $%
2^{4}=2,3^{4}=4,4^{4}=4,5^{4}=2,6^{2}=6^{4}=1$. There are not elements of
order $4$, therefore $\Theta _{7}\left( 5\right) =0$. Moreover, $\Upsilon
_{7}\left( 5\right) =0$. We obtain $\mathcal{I}_{7}\left( 5\right) =$ $%
\mathfrak{N}_{5}=2+5\ast 8=42$.

iv) If $p=11$, we have $-4=7$ and $%
2^{4}=5,3^{4}=4,4^{4}=3,5^{4}=9,6^{4}=9,7^{4}=3,8^{4}=4,9^{4}=5,10^{4}=10^{2}=1 
$. We haven't elements of order $4$, therefore $\Theta _{11}\left( 5\right)
=0$. Also, we have $\Upsilon _{11}\left( 5\right) =0$. We obtain $\mathcal{I}%
_{11}\left( 5\right) =2+\left( 11-2\right) \left( 11+1\right) =2+9\ast
12=110 $.

v) If $p=13$, we have $-4=9$ and $2^{4}=3,3^{4}=3,4^{4}=9,\mathbf{5}^{4}%
\mathbf{=1},6^{4}=9,7^{4}=9,$\newline
$\mathbf{8}^{4}\mathbf{=1},9^{4}=9,10^{4}=3,11^{4}=3,12^{4}=12^{2}=1$.
Therefore we have two elements of order $4$, namely $\{5,8\}$. It results
that $\Theta _{13}\left( 5\right) =2$ and $\Upsilon _{13}\left( 5\right) =4$%
. We obtain $\mathcal{I}_{13}\left( 5\right) =2\ast 183+5\ast \left(
50+11\ast 12\right) =1276$.

vi) If $p=17$, we have $-4=13$ and $2^{4}=16,3^{4}=13,\mathbf{4}^{4}\mathbf{%
=1},5^{4}=13,6^{4}=4,7^{4}=4,$\newline
$8^{4}=16,9^{4}=16,10^{4}=4,11^{4}=4,12^{4}=13,\mathbf{13}^{4}\mathbf{=1}%
,14^{4}=13,15^{4}=16,16^{4}=16^{2}=1$. Therefore, we have $\Theta
_{17}\left( 5\right) =2\,$, since there are only two elements of order $4$.
Also, we have $\Upsilon _{17}\left( 5\right) =4$, since the set of these
elements is $\{3,5,12,14\}$. We have $\mathcal{I}_{17}\left( 5\right) =2\ast
307+5\ast \left( 2\ast 33+15\ast 16\right) =614+1530=2144$.

vii) In the same way, for $p=19$, there are not elements of order $4$ in $%
\mathbb{Z}_{19}^{\ast }$. In this case, the number of $5$-potent elements in 
$\mathbb{H}_{\mathbb{Z}_{19}}$ is $\mathcal{I}_{19}\left( 5\right) =342$.
For $p=23$, there are not elements of order $4$ in $\mathbb{Z}_{23}^{\ast }$%
. Therefore, the number of $5$-potent elements in $\mathbb{H}_{\mathbb{Z}%
_{23}}$ is $\mathcal{I}_{23}\left( 5\right) =506$. For $p=29$, the number of 
$5$-potent elements in $\mathbb{H}_{\mathbb{Z}_{29}}$ is $\mathcal{I}%
_{29}\left( 5\right) =4872$, etc.\smallskip

We remark that from the above proposed method, we can obtain all $5$-potent
elements as the solutions of \ equations given in Case 1 and Case 2. 
\begin{equation*}
\end{equation*}%
\textbf{5.} \textbf{The number of} $k$\textbf{-potent elements in quaternion
algebra} $\mathbb{H}_{\mathbb{Z}_{p}}$%
\begin{equation*}
\end{equation*}

An element $x\in $ $\mathbb{H}_{\mathbb{Z}_{p}}$ is called $k$-potent if $k$
is the smallest positive integer such that $x^{k}=x$. From here, we have two
cases: $\mathbf{n}_{x}=0$ or $\mathbf{n}_{x}\neq 0$ and $\mathbf{n}%
_{x}^{k-1}=1$.

\textbf{Case 1,} $\mathbf{n}_{x}=0$. In this situation, we have that $%
\mathbf{t}_{x}\neq 0$, otherwise the element $x$ is nilpotent. \ From [FB;
24], Proposition 1, it results that $\mathbf{t}_{x}$ is $k$-potent in $%
\mathbb{Z}_{p}$, therefore $\mathbf{t}_{x}^{k-1}=1$ is an element of order $%
k-1$. Let $\Theta _{p}\left( k\right) $ be the number of elements of order $%
k-1$ in $\mathbb{Z}_{p}^{\ast }$.

From Proposition 1, since for a fixed $\mathbf{t}_{x}$ and $\mathbf{n}_{x}=0$%
, we have $p\left( p+1\right) $ elements in $\mathbb{H}_{\mathbb{Z}_{p}}$
satifying relation 
\begin{equation*}
\mathbf{n}_{x}=x_{0}^{2}+x_{1}^{2}+x_{2}^{2}+x_{3}^{2}=0,
\end{equation*}
it results that the number of $k$-potent elements in $\mathbb{H}_{\mathbb{Z}%
_{p}}~$with $\mathbf{n}_{x}=0,\mathbf{t}_{x}\neq 0\,$\ and $\mathbf{t}%
_{x}^{k-1}=1$ is $\Theta _{p}\left( k\right) p\left( p+1\right) $.

\textbf{Case 2,} $\mathbf{n}_{x}\neq 0$, $\mathbf{n}_{x}^{k-1}=1$ and $%
x^{k-1}=1$.

i) If $\mathbf{t}_{x}=0$, we have $x^{2}=-\mathbf{n}_{x}$ and $\mathbf{n}%
_{x}x^{k-3}+1=0$, that means $x^{k-3}=-\mathbf{n}_{x}^{-1}$.

- $k$ is even. Since $x^{2}\in \mathbb{Z}_{p}^{\ast }$, to obtain $x^{k-1}=1$%
, we must have $k-1$ even, that means $k$ is odd. Therefore, $k$ must be
odd, or $\mathbf{t}_{x}\neq 0$.

- $k$ is odd. Let $S_{k-1}=\{\mathbf{n}_{x}\in \mathbb{Z}_{p}^{\ast }$ such
that $\mathbf{n}_{x}^{k-1}=1$ and $x^{k-1}=1\}$. For each element $\alpha
\in S_{k-1}$, we must find the solutions of the equation%
\begin{equation*}
x_{1}^{2}+x_{2}^{2}+x_{3}^{2}=\alpha \text{.}
\end{equation*}

Let $\gamma _{k-1}$ be the number of elements $\alpha \in S_{k-1}$ such that 
$\alpha $ is a perfect square \ and $\delta _{k-1}$ the number of elements $%
\alpha \in S_{k-1}$ such that $\alpha $ is not a perfect square.

Therefore, in this case, we have $\Psi _{k-1}$ the number of $k$-potent
elements in $\mathbb{H}_{\mathbb{Z}_{p}}$, where: 
\begin{equation*}
\Psi _{k-1}=\left\{ 
\begin{array}{c}
\gamma _{k-1}\left( 2\mathcal{N}_{0}+\left( p-2\right) \left( p-\sin \frac{%
p\pi }{2}\right) \right) \text{, }\alpha \text{ a perfect square,} \\ 
\delta _{k-1}p\left( p-\sin \frac{p\pi }{2}\right) \text{, }\alpha \text{ is
not a perfect square}%
\end{array}%
\right. ,
\end{equation*}
with $\mathcal{N}_{0}=2p-1$\textit{, }if $p$ \textit{i}s a prime of the form 
$4l+1,l\in \mathbb{Z}$ and\textit{\ }$\mathcal{N}_{0}=1$\textit{, }if $p$ is
a prime of the form $4l+3,l\in \mathbb{Z}$.

ii) If $\mathbf{t}_{x}\neq 0$, since $x^{2}-\mathbf{t}_{x}x+\mathbf{n}_{x}=0$%
, it results $1-\mathbf{t}_{x}x^{k-2}+\mathbf{n}_{x}x^{k-3}=0$, that means 
\begin{equation*}
\mathbf{t}_{x}x^{k-2}-\mathbf{n}_{x}x^{k-3}-1=0.
\end{equation*}%
We denote $A_{k-2}=\mathbf{t}_{x}$, $B_{k-3}=-\mathbf{n}_{x},C=1$, therefore
the above relation becomes 
\begin{equation}
A_{k-2}x^{k-2}+B_{k-3}x^{k-3}+C=0.  \tag{8}
\end{equation}%
We remark that $B_{k-3}^{k-1}=\left( -1\right) ^{k-1}$\textit{mod} $p$.
Since $x^{2}=\mathbf{t}_{x}x-\mathbf{n}_{x}$, from relation $\left( 8\right) 
$, we have $A_{k-2}x^{k-4}\left( \mathbf{t}_{x}x-\mathbf{n}_{x}\right)
+B_{k-3}x^{k-3}+C=0$, then we have $\left( A_{k-2}\mathbf{t}%
_{x}+B_{k-3}\right) x^{k-3}-A_{k-2}\mathbf{n}_{x}x^{k-4}+C=0$. By using the
above notations, we have $A_{k-3}=A_{k-2}\mathbf{t}%
_{x}+B_{k-3},B_{k-4}=-A_{k-2}\mathbf{n}_{x}$. We obtain $A_{k-3}=\mathbf{t}%
_{x}^{2}-\mathbf{n}_{x},B_{k-3}=-\mathbf{t}_{x}\mathbf{n}_{x}$. Step by
step, we have $A_{2}x^{2}+B_{1}x+C=0$. It results the system 
\begin{equation*}
\left\{ 
\begin{array}{c}
x^{2}-\mathbf{t}_{x}x+\mathbf{n}_{x}=0 \\ 
A_{2}x^{2}+B_{1}x+C=0%
\end{array}%
\right.
\end{equation*}
and we have $\left( A_{2}\mathbf{t}_{x}+B_{1}\right) x=A_{2}\mathbf{n}_{x}+C$%
. If $A_{2}\mathbf{t}_{x}+B_{1}\neq 0$, we have $x\in \mathbb{Z}_{p}$ and it
is an element of order $k-1$. If \ $A_{2}\mathbf{t}_{x}+B_{1}=0$, then $A_{2}%
\mathbf{n}_{x}+1=0$, since $C=1$. It is clear that $A_{2}\neq 0$, therefore $%
\mathbf{n}_{x}=-A_{2}^{-1}$, then $A_{2}^{k-1}=1$. Also, $\mathbf{t}%
_{x}=-B_{1}A_{2}^{-1}=B_{1}\mathbf{n}_{x}$, then $\mathbf{t}%
_{x}^{k-1}=B_{1}^{k-1}$ and $\mathbf{n}_{x}=B_{1}^{-1}\mathbf{t}_{x}$. For $%
x=x_{0}+x_{1}\mathbf{i}+x_{2}\mathbf{j}+x_{3}\mathbf{k}$, we have $%
x_{0}=2^{-1}\mathbf{t}_{x}$ and relation $\mathbf{n}_{x}=B_{1}^{-1}\mathbf{t}%
_{x}$, with $\mathbf{t}_{x}^{k-1}=B_{1}^{k-1}$, implies 
\begin{equation*}
\left( 2^{-1}\mathbf{t}_{x}\right)
^{2}+x_{1}^{2}+x_{2}^{2}+x_{3}^{2}=B_{1}^{-1}\mathbf{t}_{x}\text{,}
\end{equation*}%
with $\left( B_{1}^{-1}\mathbf{t}_{x}\right) ^{k-1}=1$.

We obtain the equation 
\begin{equation}
x_{1}^{2}+x_{2}^{2}+x_{3}^{2}=B_{1}^{-1}\mathbf{t}_{x}-4^{-1}\mathbf{t}%
_{x}^{2}\text{.}  \tag{9}
\end{equation}%
Let $\mathcal{N}_{k-1}$ be the number of solutions of equation $\left(
9\right) $. Therefore the number $k$-potent elements is $\mathcal{I}%
_{p}\left( k\right) =\Theta _{p}\left( k\right) p\left( p+1\right) +\Psi
_{k-1}+\mathcal{N}_{k-1}$.\smallskip

\textbf{Theorem 21.} \textit{With the above notations, we have} 
\begin{equation}
\mathcal{I}_{p}\left( k\right) =\Theta _{p}\left( k\right) p\left(
p+1\right) +\Psi _{k-1}+\mathcal{N}_{k-1}.  \tag{10}
\end{equation}%
We remark that from the above proposed method, we can obtain all $k$-potent
elements as the solutions of \ equations given in Case 1 and Case 2.%
\begin{equation*}
\end{equation*}%
\textbf{6.} \textbf{An application}%
\begin{equation*}
\end{equation*}

In the paper [Ni; 42], the author studied the existence of the $n$-roots of
a real quaternion, namely, was found the solutions and number of the
solutions of the equation 
\begin{equation}
x^{n}=a  \tag{11}
\end{equation}%
in real quaternions. If $a$ is not a real number, the above equation has
exactly $n$ distinct roots. If $\alpha $ is a real number and $n\neq 2$,
equation has infinitely many solutions. If $n=2,\alpha >0$, we have oly two
square roots: $\sqrt{\alpha }$ and $-\sqrt{\alpha }$.

In the following, we will study the behaviour of such an equation over
finite fields. In the following, we will study and count the number of
solutions of equation 
\begin{equation}
x^{k}=1  \tag{12}
\end{equation}%
in $\mathbb{H}_{\mathbb{Z}_{p}}$. We consider $\ k$ a positive integer and $%
1<d_{1}<d_{2}<...<d_{r-1}<d_{r}=k$ the divisors of $k$. We remark that $1$
is a solution of equation $\left( 12\right) $. Also, we remark that the
solutions of equations 
\begin{equation}
x^{d_{j}}=1,j\in \{1,2,...r\}  \tag{13}
\end{equation}%
are also solutions for equation $\left( 12\right) $. We remark that the
solutions of equation $\left( 13\right) $ are all $(d_{j}+1)$-potent
elements with norm different from zero. Let $\mathcal{S}_{d_{j}}$ be the
number of solutions of the equation $\left( 13\right) $ and $\mathcal{N}%
_{p}\left( k\right) $ the number of all solutions of equation $\left(
12\right) $.

From Theorem 21, we obtain that $\mathcal{N}_{p}\left( k\right) =\mathcal{S}%
_{1}+\underset{j\in \{1,2,...r\}}{\sum }\mathcal{S}_{d_{j}}$. Since $%
\mathcal{S}_{1}=1$ and $\mathcal{S}_{d_{j}}=\mathcal{I}_{p}\left(
d_{j}+1\right) -\Theta _{p}\left( d_{j}+1\right) p\left( p+1\right) $, we
have $\mathcal{N}_{p}\left( k\right) =\mathcal{S}_{1}+\underset{j\in
\{1,2,...r\}}{\sum }\mathcal{S}_{d_{j}}=$\newline
$=1+\underset{j\in \{1,2,...r\}}{\sum }\left( \mathcal{I}_{p}\left(
d_{j}+1\right) -\Theta _{p}\left( d_{j}+1\right) p\left( p+1\right) \right) $%
. Therefore, we obtain the following theorem.\smallskip

\textbf{Theorem 22.} \textit{The number of solutions of equation} $\left(
12\right) $ \textit{is} 
\begin{equation}
\mathcal{N}_{p}\left( k\right) =1+\underset{j\in \{1,2,...r\}}{\sum }\left( 
\mathcal{I}_{p}\left( d_{j}+1\right) -\Theta _{p}\left( d_{j}+1\right)
p\left( p+1\right) \right) .  \tag{14}
\end{equation}

\textbf{Example 23.} \thinspace i) For $p=3$, $k=4$, we count $\mathcal{N}%
_{3}\left( 4\right) $, the number of solutions of equation 
\begin{equation*}
x^{4}=1\text{. }
\end{equation*}%
We have $4=2^{2}$, therefore $1<2<4$ are the divisors of $4$. We count the $%
3 $-potent and $5$-potent elements in $\mathbb{H}_{\mathbb{Z}_{3}}$ with
nonzero norm. We use Theorem 14 and Theorem 19. It results $\mathcal{N}%
_{3}\left( 4\right) =1+\left( \mathcal{I}_{3}\left( 3\right) -\Theta
_{3}\left( 3\right) \ast 3\ast 4\right) +(\mathcal{I}_{3}\left( 5\right)
-\Theta _{3}\left( 5\right) \ast 3\ast 4)=$\newline
$=1+(3^{2}+3+1)+(2+1\ast 4)=20$.

ii) For $p=3$, $k=6$, we count $\mathcal{N}_{3}\left( 6\right) $, the number
of solutions of equation 
\begin{equation*}
x^{6}=1\text{. }
\end{equation*}%
We have $6=2\ast 3$, therefore $1<2<3<6$ are the divisors of $6$. We count
the $3$-potent, $4$-potent and $7$-potent elements in $\mathbb{H}_{\mathbb{Z}%
_{3}}$ with nonzero norm. We use Theorem 14 and Example 18 i). We have $%
(3^{2}+3+1)=13$ elements which are $3$-potents with nonzero norm and $8$
elements which are $4$-potents with nonzero norm. By using the similar
procedure, we count the $7$-potents with nonzero norm. This number is $8$.
Therefore, we have $\mathcal{N}_{3}\left( 6\right) =1+13+8+8=30$.

iii) For $p=5$, $k=4$, we count $\mathcal{N}_{5}\left( 4\right) $, the
number of solutions of equation 
\begin{equation*}
x^{4}=1\text{. }
\end{equation*}%
As above, we count the $3$-potent and the $5$-potent elements in $\mathbb{H}%
_{\mathbb{Z}_{5}}$ with nonzero norm. We use Theorem 14 and Theorem 19. It
results $\mathcal{N}_{5}\left( 4\right) =1+$ $(5^{2}+5+1)+$ $[2+5\ast \left(
2\ast 9+3\ast 4\right) ]=184$.

\begin{equation*}
\end{equation*}

\textbf{Conclusions.} \ In this paper, we presented a descriptive formula to
count the number of all of $k$-potent elements over $\mathbb{H}_{\mathbb{Z}%
_{p}}$ and, for $k\in \{3,4,5\}$, we give an explicit formula for these
values. The proposed methods give us all elements which are tripotents, $4$%
-potents, $5$-potents or, in general, $k$-potents, $k$ a positive integer,
\thinspace $k\geq 6$. Moreover, as an application, we count the number of
solutions of the equation $x^{k}=1$ over $\mathbb{H}_{\mathbb{Z}_{p}}$ and
we gave a lot of examples. For this purpose, we used computer as a tool to
check and understand the behavior of these elements in each studied case.
Due to a lot of computation involved, this allowed us to give a correct
mathematical proof for obtained relations. By using C + + software (https:

//www.programiz.com/cpp-programming/online-compiler/, accessed on 20
November 2024), we obtain, for example, the below code:\medskip 

\textbf{C++ Code: Counting 3-Potent Quaternions}
\begin{verbatim}
	// Iterate over all quaternions in Z7
	for (int a = 0; a < 7; a++) {
		for (int b = 0; b < 7; b++) {
			for (int c = 0; c < 7; c++) {
				for (int d = 0; d < 7; d++) {
					int q[4] = {a, b, c, d};
					int q2[4], q3[4];
					
					quaternion_power(q, 2, q2); // Compute Q^2
					quaternion_power(q, 3, q3); // Compute Q^3
					
					// Check if Q^3 = Q and Q^2 != Q
					if (q3[0] == q[0] && q3[1] == q[1] && q3[2] == q[2] && q3[3] == q[3]) {
						if (q2[0] != q[0] || q2[1] != q[1] || q2[2] != q[2] || q2[3] != q[3]) {
							count++;
						}
					}
				}
			}
		}
	}
	
	std::cout << "Number of quaternions that are 3-potent but not idempotent in Z7: "
	<< count << std::endl;	
	113
\end{verbatim}

\textbf{C++ Code: Counting 4-Potent Quaternions}
\begin{verbatim}
	int count = 0;
	// Iterate over all possible quaternions in Z11
	for (int a = 0; a < 11; a++) {
		for (int b = 0; b < 11; b++) {
			for (int c = 0; c < 11; c++) {
				for (int d = 0; d < 11; d++) {
					int q[4] = {a, b, c, d};
					int q2[4], q3[4], q4[4];
					
					quaternion_power(q, 2, q2); // Compute Q^2
					quaternion_power(q, 3, q3); // Compute Q^3
					quaternion_power(q, 4, q4); // Compute Q^4
					
					// Check if Q^4 = Q but Q^2 != Q and Q^3 != Q
					if (q4[0] == q[0] && q4[1] == q[1] && q4[2] == q[2] && q4[3] == q[3]) {
						if ((q2[0] != q[0] || q2[1] != q[1] || q2[2] != q[2] || q2[3] != q[3]) &&
						(q3[0] != q[0] || q3[1] != q[1] || q3[2] != q[2] || q3[3] != q[3])) {
							count++;
						}
					}
				}
			}
		}
	}
	
	std::cout << "Number of quaternions that are 4-potent but not idempotent or tripotent in Z11: "
	<< count << std::endl;
	110
\end{verbatim}

\textbf{C++ Code: Counting 5-Potent Quaternions}
\begin{verbatim}
 // Function to calculate the total number of quaternions satisfying Q^5 = Q
	int count_quaternions_Q5(int mod) {
		int count = 0;  // Initialize a counter for valid solutions
  // Iterate over all quaternions in Z13
		// Iterate through all possible values of the quaternion components (a, b, c, d) in Zmod
		for (int a = 0; a < mod; ++a) {
			for (int b = 0; b < mod; ++b) {
				for (int c = 0; c < mod; ++c) {
					for (int d = 0; d < mod; ++d) {
						// Construct the quaternion Q = (a, b, c, d)
						tuple<int, int, int, int> q = make_tuple(a, b, c, d);
						
						// Calculate Q^5 using the quaternion_power function
						auto q5 = quaternion_power(q, 5, mod);
						
						// Check if Q^5 is equal to Q
						if (q5 == q) {
							++count;  // If true, increment the counter
						}
					}
				}
			}
		}
		
		// Return the total number of solutions found
		return count;
	}
	1276
\end{verbatim}

As a further research, we intend to rafinate the formula for $k$-potent
quaternions and to obtain exact formulae for other values of $k$. 
\begin{equation*}
\end{equation*}

\textbf{References}%
\begin{equation*}
\end{equation*}

[AD; 12] Aristidou M., Demetre A., \textit{A Note on Nilpotent Elements in
Quaternion Rings over } $\mathbb{Z}_{p}$, International Journal of Algebra,
6(14)(2012), 663 - 666.

[FH; 58] Fine N. J., Herstein I. N., \textit{The Probability that a Matrix
be Nilpotent}, Illinois Journal of Mathematics, 2(4A) (1958), 499-504.

[FB; 24] Flaut C., Baias A., \textit{Some Remarks Regarding Special Elements
in Algebras Obtained by the Cayley--Dickson Process over} $\mathbb{Z}_{p}$,
Axioms 3(6)(2024), 351; https://doi.org/10.3390/axioms13060351

\textbf{\ [}FSF; 19\textbf{] }Flaut C, Ho\v{s}kov\'{a}-Mayerov\'{a} \v{S},
Flaut D, \textit{Models and Theories in Social Systems}, Springer Nature,
2019, https://www.springer.com/us/book/9783030000837, e-ISBN:
978-3-030-00084-4, ISBN: 978-3-030-00083-7, DOI: 10.1007/978-3-030-00084-4

[MS; 11] Miguel C. J., Ser\^{o}dio R., \textit{On the Structure of
Quaternion Rings over} $\mathbb{Z}_{p}$, International Journal of Algebra,
5(27)(2011), 1313 - 1325.

[Ni; 42] Niven I., \textit{The roots of a quaternion}, The American
Mathematical Mounthly, 49(6)(1942), 386-388.

\textbf{[}Sc; 66\textbf{]} Schafer R. D., \textit{An Introduction to
Nonassociative Algebras,} Academic Press, New-York, 1966.

\textbf{[}Vo; 21\textbf{]} Voight J, \textit{Quaternion Algebras}, Springer
Nature Switzerland AG, ISBN 978-3-030-56692-0

[W; 22] Wongkumpra P., \textit{Solutions to quadratic equations over finite
fields}, Master of Science thesis, 2022,
https://digital\_collect.lib.buu.ac.th/dcms/files/62910234.pdf%
\begin{equation*}
\end{equation*}%
Cristina FLAUT

Faculty of Mathematics and Computer Science, Ovidius University,

Bd. Mamaia 124, 900527, Constan\c{t}a, Rom\^{a}nia,

http://www.univ-ovidius.ro/math/

e-mail: cflaut@univ-ovidius.ro; cristina\_flaut@yahoo.com

\begin{equation*}
\end{equation*}%
Andreea BAIAS

PhD student at Doctoral School of Mathematics,

Ovidius University of Constan\c{t}a, Rom\^{a}nia,

e-mail: andreeatugui@yahoo.com

\end{document}